\documentclass[11pt]{amsart}

\def\cal{\mathcal}
\def\F{\mathcal{F}}
\def\E{\mathcal{E}}
\def\H{{H}}
\def\R{\mathcal{R}}
\def\B{\mathcal{B}}
\def\A{\mathcal{A}}
\def\la{\langle}
\def\ra{\rangle}
\def\u{\underline}
\def\sf{Sub^\F}
\def\pf{Pos^\F}
\def\at{\alpha_t}
\def\La{\Leftarrow}
\def\laa{\leftarrow}
\def\Ra{\Rightarrow}
\def\P1{\mathcal{P}}
\newcommand{\natur}{\mbox{\it{I\hspace{-0.2em}N}}}
\def\nn{\natur}

\begin{document}

\email{shtrakov@aix.swu.bg}
\address{Department of Computer Science \\
South-West University,  2700 Blagoevgrad, Bulgaria}
\title{Multi-Solid Varieties and Mh-transducers}%Paper title

\author[Sl. Shtrakov]{Slavcho Shtrakov}
\urladdr{http://home.swu.bg/shtrakov}
\date{}

\date{}

\theoremstyle{plain}
\newtheorem{theorem}{Theorem}
\newtheorem{lemma}{Lemma}
\newtheorem{proposition}{Proposition}
\newtheorem{corollary}{Corollary}
\newtheorem{definition}{Definition}
\theoremstyle{definition}
\newtheorem{example}{Example}
\newtheorem{remark}{Remark}

\maketitle

\begin{abstract}
We consider the concepts of colored terms and
multi- hypersubstitutions. If $t\in W_\tau(X)$ is a term of type
$\tau$, then any mapping  $\alpha_t:Pos^\F(t)\to \nn$ of the
non-variable positions of a term into the set of natural numbers is
called a coloration  of  $t.$ The set $W_\tau^c(X)$ of  colored
terms consists of all pairs $\langle t,\alpha_t\rangle.$
Hypersubstitutions are maps which assign to each operation symbol a
term with the same arity.
 If
$M$ is a monoid of hypersubstitutions then any sequence $\rho =
(\sigma_1,\sigma_2,\ldots)$ is   a mapping $\rho:\nn\to M$, called a
multi-hypersubstitution over $M$.  An identity $t\approx s$,
satisfied in a variety $V$ is an $M$-multi-hyperidentity if its
images $\rho[t\approx s]$ are  also satisfied in $V$ for all
$\rho\in M$. A variety $V$ is $M$-multi-solid, if all its identities
are $M-$multi-hyperidentities. We prove a series of inclusions and
equations concerning $M$-multi-solid varieties.  Finally we give an
automata realization of multi-hypersubstitutions and colored terms.
\end{abstract}

\section*{Introduction}

Let ${\F}$ be a set of  {\it operation\ symbols}, and
$\tau:{{\F}}\to N$ be a  {\it type} or {\it signature}.

 Let $X$  be a finite set of variables, then the
set $W_{\tau}(X)$ of {\it terms\ of\ type\ $\tau$ } with\ variables\
from\ $X $ is the smallest set, such that
\\ \hspace*{1cm}
$(i)$\quad $X\cup\F_0\subseteq W_\tau(X);$
\\ \hspace*{1cm}
$(ii)$\quad if $f $ is $n-$ary operation symbol and
$t_1,\ldots,t_{n}$ are terms, then the ``string" $f(t_1\ldots
t_{n})$ is a term.

 An algebra
 ${\cal A} = \la A; \F^{\A}\ra$
  of type $\tau$ is a pair consisting of a set $A$ and a
   set $\F^\A$ of operations defined on $A$.
  If $f\in\F$, then $f^\A$ denotes a $\tau(f)$-ary operation on the
set $A.$
     An identity $s \approx t$ is satisfied
                in the algebra ${\cal A}$    (written $\A\models s \approx
                t$), if
                $s^\A= t^\A.$
          $\A lg(\tau)$ denotes the class of all algebras
   of type $\tau$ and  $Id(\tau)$ - the set of all identities
of type $\tau.$
     The pair $(Id, Mod)$
     is a Galois connection  between the classes of algebras from
       $\A lg(\tau)$ and subsets of $Id(\tau)$,
where $Id(\R):=\{t\approx s\ | \forall
\A\in\R, \ (\A\models t\approx s)\}$ and $Mod(\Sigma):=\{\A |
\forall t\approx s\in\Sigma, \ (\A\models t\approx s)\}$.
 The fixed points
 with respect to the closure operators
$Id Mod$ and $Mod Id$ form complete lattices
\\
\centerline{ ${\cal L}(\tau):= \{\R \mid \R \subseteq \A lg(\tau)
~\mbox{and}~ ModId\R =\R\} \mbox{  and}$}\\
\centerline{ ${\cal E}(\tau):= \{\Sigma \mid \Sigma
 \subseteq Id(\tau) ~\mbox{and}~ IdMod\Sigma = \Sigma\}$}
of all varieties of type $\tau$ and of all equational theories
(logics) of type $\tau$. These lattices
 are dually isomorphic.

 A {\it hypersubstitution of type }$\tau$
(briefly {\it  a hypersubstitution}) is a mapping which assigns to
each operation symbol $f\in\F$ a term   $\sigma(f)$ of type $\tau$,
which has the same arity as the operation symbol  $f$ (see
\cite{den51}). The set of all hypersubstitutions of type  $\tau$ is
denoted by $Hyp(\tau).$
 If  $\sigma$ is a hypersubstitution, then it can be uniquely
  extended to a mapping
  $\hat\sigma:W_\tau(X)\to W_{\tau}(X)$
  on the set of all
  terms of type $\tau$, as follows

\begin{tabular}{ll}
$(i)$&if $t = x_j$ \  for some  $j \geq 1$, \
 then $\hat\sigma[t] = x_j$; \\
$(ii)$& if $t= f(t_1,\ldots,t_{n})$ ,\  then $\hat\sigma[t] =
\sigma(f)(\hat\sigma[t_1],\ldots, \hat\sigma[t_{n}])$,
\end{tabular}
\\
where $f$ is an  $n$-ary operation symbol and
  $t_1$,$\ldots,t_{n}$ are terms.

The set  $Hyp(\tau)$ is a monoid.
%%%%%%%%%%%%%%%%%%%%%%%%%

 Let ${  M} $ be a submonoid of $Hyp(\tau)$.
An algebra ${\cal A}$ is said to $M$-hypersatisfy
 an identity $t \approx s$ if for every
hypersubstitution $\sigma \in M$, the identity ${\hat
\sigma}[t]\approx{\hat \sigma}[s]$ holds in ${\cal A}$. A variety
$V$ is called {\it $M$-solid} if every identity of $V$ is
$M$-hypersatisfied in  $V$.

The closure operator is defined on the  set of identities of a given
type $\tau$  as follows: $\chi_M[u \approx v]:= \{\hat\sigma[u]
\approx \hat\sigma[v] \mid \sigma \in M\}$\ \mbox{ and}\
$\chi_M[\Sigma] = \bigcup \limits_{u \approx v \in \Sigma}\chi_M[u
\approx v].$

Given an algebra ${\cal A} = \la A; \F^{\cal A}\ra$ and a
hypersubstitution $\sigma$, then  $\sigma[{\cal A}] = \la A;
({\sigma(\F^{\cal A})}\ra:= \la A; (\sigma(f)^{\cal A})_{f\in
\F}\ra$ is called the {\it derived algebra}. The closure operator
$\psi_M$ on the set of algebras of a given type $\tau$, is defined
as follows: $\psi_M[{\cal A}] = \{\sigma[{\cal A}] \mid \sigma \in M
\}$\ \mbox{ and}\  $\psi_M[\R] = \bigcup
\limits_{{\cal A} \in \R}\psi_M[{\cal A}].$ \\

It is well known \cite{den51} that if  ${ M}$ is a monoid of
hypersubstitutions of type $\tau$, then  the class
 of all $M$-solid varieties of type $\tau$ forms
a complete sublattice of the lattice ${\cal L}(\tau)$ of all
varieties of type $\tau$.

 Our aim is to transfer these results to another kind of
hypersubstitution and to coloured terms.

In Section \ref{sec2} we present two constructions which  produce
composed terms. The first one is inductive and the resulting term is
obtained by
 simultaneous replacement of a
subterm in all places where it occurs in a given term with another
term of the same type. The positional composition gives a composed
term as a result of the replacement of subterms in given positions
with other terms of the same type.  The positional composition of
colored terms is an associative operation. In \cite{dks1} the
authors studied colored terms which are supplied with one
coloration. This is a very ``static" concept where each term has one
fixed coloration. Here we consider composition of terms which
produces an image of  terms and coloration of this image. This
``dynamical" point of view gives us an advantage when studying
multi-hypersubstitutions, multi-solid varieties etc.

In Section \ref{sec4} we use colored terms to investigate the monoid
of multi-hypersubstitutions. It is proved that the lattice of
multi-solid varieties is a sublattice of the lattice of solid
varieties. A  series of
 assertions are proved, which characterize  multi-solid varieties and the
corresponding closure operators. We study  multi-solid varieties
by deduction of a fully invariant congruence. The completeness
theorem for multi-hyperequational theories is proved.

A tree automata realization of multi-hypersubstitutions is given in
 Section \ref{sec7}.

\section{Composition  of colored
terms}\label{sec2}

The concept of the composition of mappings is  fundamental in
almost all mathematical theories. Usually we consider composition
as an operation which inductively replaces some variables with
other objects such as functions, terms, etc. Here, we consider a
more general case when the replacement can be applied to objects
which may be variables or subfunctions, subterms,
 etc. which are located at a  given set of positions.

If $t$ is a term, then the
set $var(t)$ consisting of these elements of $X$ which occur in
$t$ is called the set of
 {\it input variables (or variables)} for $t$.  If
$t=f(t_1,\ldots,t_n)$ is a non-variable term, then $f$ is {\it
root symbol (root)} of $t$ and we will write $f=root(t).$
%\begin{definition}%
\label{d-Sub} For a term  $t\in W_\tau(X)$ the set   $Sub(t)$ of
its subterms
 is defined as follows:
if $t\in X\cup \F_0$, then $Sub(t)=\{t\}$ and if
$t=f(t_1,\ldots,t_{n})$, then $Sub(t)=\{t\}\cup
Sub(t_1)\cup\ldots\cup Sub(t_n).$
%\end{definition}

The $depth$ of a term $t$ is defined inductively: if $t\in
X\cup\F_0$ then $Depth(t)=0;$ \ and  if $t=f(t_1,\ldots,t_{n})$,
then\\
 $Depth(t)=max \{Depth(t_1),\ldots,Depth(t_{n})\} +1.$
\begin{definition}%
\label{d-indkompGeneral} Let $r,s,t\in W_\tau(X)$ be three terms
of type $\tau$. By $t(r \laa s)$ we will denote the term, obtained
by simultaneous replacement of  every occurrence of $r $ as a
 subterm of $t$ by $s$.  This term is called
 the
{\it inductive composition } of the terms $t$ and $s $, by $r $.
I.e.
\\
\begin{tabular}{ll}
 \hspace*{.5cm}&
 $(i)$\quad $t(r\laa s)=t$ if $r\notin
Sub(t);$
\\ \hspace*{.5cm}&
 $(ii)$\quad $t(r\laa s)=s$ if $t=r$  and
\\ \hspace*{.5cm}&
$(iii)$\quad $t(r\laa s)=f(t_1(r\laa s),\ldots,t_n(r\laa s)),$
\\ \hspace*{.5cm}&\hspace*{.7cm} if
$t=f(t_1,\ldots,t_n)$ and $r\in Sub(t)$, $r\neq t$.
\end{tabular}
\end{definition}
If $r_i\notin Sub(r_j)$ when $i\neq j$, then $t(r_1\laa
s_1,\ldots,r_m\laa s_m)$ means the inductive composition of
$t,r_1,\ldots,r_m,s_1,\ldots,s_m$. In the particular case when
$r_j=x_j$ for $j=1,\ldots,m$ and $var(t)=\{x_1,\ldots,x_m\}$ we
will briefly write $t(s_1,\ldots,s_m)$ instead of $t(x_1\laa
s_1,\ldots,x_m\laa s_m)$.

 Let $\tau$ be a type and $\F$ be its set of operation
symbols. Denote by $maxar=\max\{\tau(f)|f\in\F\}$ and
$\nn_\F:=\{m\in\nn\ |\ m\leq maxar\}$.  Let $\natur_\F^*$ be the
set of all finite strings over  $\natur_\F.$ The set $\natur_\F^*$
is naturally ordered by $p\preceq q \iff p$\ is a prefix of $q.$
The Greek letter $\varepsilon$, as usual  denotes the empty word
(string) over $\natur_\F.$

For any term $t$, the set of  positions $Pos(t)\subseteq
\natur_{\F}^*$ of $t$ is inductively defined as follows:
 $Pos(t)=\{\varepsilon\}$ if $t\in X\cup \F_0$ and
$Pos(t):=\{\varepsilon\}\bigcup_{1\leq i\leq n}(i Pos(t_i)) ,$ \ \
if $t=f(t_1,\ldots,t_n),$ $n\geq0$, where $i Pos(t_i):=\{i q |
q\in Pos(t_i)\},$ and $iq$ is concatenation of the strings $i$ and
$q$ from
 $\natur_\F^*.$

For a given position $p\in Pos(t)$,  the length of $p$ is denoted
by  $l(p)$.

Any term can be regarded  as a tree with nodes labelled with the
operation symbols and its leaves   labelled as variables or nullary
operation symbols.

  Let $t\in W_\tau(X)$ be a term of type $\tau$ and let $sub_t:Pos(t)\to
Sub(t)$  be the function  which maps each position in a term $t$ to
the subterm of $t$, whose root node occurs at that position.

\begin{definition}%
\label{d-kompos3}
 Let $t,r\in W_\tau(X)$ be two terms of type
$\tau$ and $p\in Pos(t)$ be a position in $t.$  The positional
composition of $t$ and $r$ on $p$ is a term $s:=t(p;r)$ obtained
from $t$ when replacing   the term $sub_t(p)$ by $r$  on the
position $p$, only.
\end{definition}
More generally, the positional composition of terms is naturally
defined for the compositional pairs
$t(p_1,\ldots,p_m;t_1,\ldots,t_m)$, also.

\begin{remark}\label{rem-assoc}%\cite{shau}
The positional composition has the following properties:

 1.
If $\la\la p_1,p_2\ra,\la t_1,t_2\ra\ra$ is a compositional pair
of $t$, then
$$t(p_1,p_2;t_1,t_2)=t(p_1;t_1)(p_2;t_2)=t(p_2;t_2)(p_1;t_1);$$

 2. If $S=\la p_1,\ldots,p_m\ra$ and $T=\la t_1,\ldots,t_m\ra$
with  $$(\forall
 p_i,p_j\in S)\    (i\neq j \implies p_i\not\prec p_j\  \&\  p_j\not\prec p_i)$$
 and $\pi$ is a permutation of the set $\{1,\ldots,m\}$, then
$$t(p_1,\ldots,p_m;t_1,\ldots,t_m)=
t(p_{\pi(1)},\ldots,p_{\pi(m)};t_{\pi(1)},\ldots,t_{\pi(m)}).$$

 3.If $t,s,r\in W_\tau(X)$, \ $p\in Pos(t)$ and $q\in
Pos(s)$, then $t(p;s(q;r))=t(p;s)(pq;r)$.
\end{remark}

We will denote respectively  $``\laa"$ for  inductive composition
and $``;"$ for  positional composition.

%\label{sec3}
The concept of colored terms is important when studying deductive
closure of sets of identities and the lattice of varieties of a
given type. Colored terms (trees) are useful tools   in Computer
Science, General Algebra, Theory of Formal Languages, Programming,
Automata theory etc.

Let $t$ be a term of type $\tau.$ Let us denote by $\sf(t)$ the
set of all subterms of $t$ which are not variables, i.e. whose
roots are labelled by an operation symbol from $\F$ and let
$\pf(t):=\{p\in Pos(t)\ |\ sub_t(p)\in\sf(t)\}$. Any function
$\at:\pf(t)\to \nn$ is called a {\it  coloration} of the term $t.$

For a given term $t$,   $Pos^X(t)$ denotes the set of all its
variable  positions i.e. $Pos^X(t):=Pos(t)\setminus\pf(t).$

By $C_t$ we denote the set of all colorations of the term $t$ i.e.
$C_t:=\{\at\ |\ \at:\pf(t)\to \nn\}$. If $p\in\pf(t)$ then
$\alpha_t(p)\in \nn$ denotes the value of the function $\alpha_t$
which is associated with the root operation symbol of the subterm
$s=sub_t(p)$, and $\alpha_t[p]\in C_s$ denotes the "restriction"
of the function $\alpha_t$ on the  set $\pf(s)$ defined by
$\alpha_t[p](q)=\at(pq)$ for all $q\in\pf(s)$.

\begin{definition}\label{d-colorTerms} The set $W^c_\tau(X)$ of
all colored terms of type $\tau$ is defined as follows:

$(i)$\quad $X\subset W^c_\tau(X)$;

$(ii)$\quad If $f\in\F$, then $\la f,q\ra\in W^c_\tau(X)$ for each
$q\in \nn;$

$(iii)$\quad If $t=f(t_1,\ldots,t_n)\in W_\tau(X)$, then $\la
t,\at\ra\in W^c_\tau(X)$ for each $\at\in C_t.$
\end{definition}\

Let $\la t,\at\ra\in W^c_\tau(X)$. The set $Sub_c(\la t,\at\ra)$
of colored subterms of $\la t,\at\ra$ is defined as follows:

 For
$x\in X$ we have $Sub_c(x):=\{x\}$ \\
and if  $\la t,\at\ra=\la
f,q\ra(\la t_1,\alpha_{t_1}\ra,\ldots,\la t_n,\alpha_{t_n}\ra)$,
then
$$Sub_c(\la
t,\at\ra):=\{\la t,\at\ra\}\cup Sub_c(\la
t_1,\alpha_{t_1}\ra)\cup\ldots\cup Sub_c(\la
t_n,\alpha_{t_n}\ra).$$

 Let
 $\la t,\alpha_t\ra,$
$\la r,\alpha_{r}\ra$ and $\la s,\alpha_{s}\ra$ be colored terms
of type $\tau$. Their inductive composition $\la t,\alpha_t\ra(\la
r,\alpha_{r}\ra\laa \la s,\alpha_{s}\ra)$ is defined as follows:

 $(i)$   if $t=x_i\in X$, then
 $$x_i(\la
r,\alpha_{r}\ra\laa \la s,\alpha_{s}\ra):=\left\{\begin{array}{ll}
                \la s,\alpha_{s}\ra & \mbox{if}\ \   r=x_i;\\
                x_i & \mbox{otherwise }\

\end{array}
\right.$$

  $(ii)$  if  $t=f(t_1,t_2,\ldots,t_n)$, $\at(\varepsilon)=q,$
$\at[i]=\alpha_{t_i}$  for $i=1,2,\ldots,n$,
  then
    $\la t,\alpha_t\ra(\la r,\alpha_{r}\ra\laa \la
s,\alpha_{s}\ra):=
                \la s,\alpha_{s}\ra $ when $ \la r,\alpha_{r}\ra=\la t,\alpha_t\ra$
                and
                \\ \centerline{$\la t,\alpha_t\ra(\la r,\alpha_{r}\ra\laa \la
s,\alpha_{s}\ra):=$}
\\ \centerline{$ \la f,q\ra(\la t_1,\alpha_{t_1}\ra(\la
r,\alpha_{r}\ra\laa \la s,\alpha_{s}\ra),\ldots, \la
t_n,\alpha_{t_n}\ra(\la r,\alpha_{r}\ra\laa \la s,\alpha_{s}\ra)$}
otherwise.

If $\la r_i,\alpha_{r_i}\ra\notin Sub_c(\la r_j,\alpha_{r_j}\ra)$
for $i\neq j$, then the denotations \\
 $\la t,\alpha_t\ra(\la
r_1,\alpha_{r_1}\ra\laa \la s_1,\alpha_{s_1}\ra,\ldots, \la
r_m,\alpha_{r_m}\ra\laa \la s_m,\alpha_{s_m}\ra)$ is clear.

Let $\la t,\at\ra,\la s,\alpha_s\ra\in W^c_\tau(X)$ be two colored
terms of type $\tau$ and let $ p\in Pos(t)$.  The positional
composition
 of the colored terms $\la t,\at\ra$ and $\la
s,\alpha_{s}\ra $ at the position $p$ is defined as follows:
 \\
\centerline{ $\la t,\at\ra(p;\la s,\alpha_{s}\ra):=\la t(p;s),\alpha\ra$,}\\
where
\[ \alpha(q)=\left\{\begin{array}{ll}
                \alpha_s(k) & \mbox{if   }\ \
                   q=pk,\ \mbox{for some } k\in Pos(s) ;\\
                %&\\
                \alpha_{t}(q)& \mbox{otherwise. }\
                \end{array}
\right.\]

The positional composition of colored terms can be defined for a
 sequence  $(p_1,\ldots,p_m)\in Pos(t)^m$ of positions
with  $$(\forall
 p_i,p_j\in S)\    (i\neq j \implies p_i\not\prec p_j\  \&\  p_j\not\prec p_i).$$ It is denoted by\\
\centerline{ $\la t,\alpha_t\ra(p_1,\ldots,p_m;\la
s_1,\alpha_{s_1}\ra,\ldots,\la s_m,\alpha_{s_m}\ra)$.}

\begin{theorem}
\label{t-ColorCompos1}  If $t,s,r\in W_\tau(X)$, \ $p\in Pos(t)$ and
$q\in Pos(s)$, then $$\la t,\alpha_t\ra(p;\la s,\alpha_s\ra(q;\la
r,\alpha_r\ra))=\la t,\alpha_t\ra(p;\la s,\alpha_s\ra)(pq;\la
r,\alpha_r\ra),$$ where $\alpha_t\in C_t$, $\alpha_s\in C_s$,
$\alpha_r\in C_r$.
\end{theorem}
\begin{proof}
Let us consider the non-trivial case when $t,r$ and $s$ are not
variables. Thus we obtain $\la s,\alpha_s\ra(q;\la
r,\alpha_r\ra)=\la s(q;r),\alpha\ra$, where
\[ \alpha(m)=\left\{\begin{array}{ll}
                \alpha_r(k) & \mbox{if   }\ \
                   m=qk,\ \mbox{for some } k\in \pf(r) ;\\
                %&\\
                \alpha_{s}(m)& \mbox{otherwise}\
                \end{array}
\right.\] and $\la t,\at\ra(p;\la s(q;r),\alpha\ra)=\la
t(p;s(q;r)),\beta\ra$, where
\[ \beta(l)=\left\{\begin{array}{ll}
                \alpha(v) & \mbox{if   }\ \
                   l=pv,\ \mbox{for some } v\in \pf(s(q;r)) ;\\
                %&\\
                \alpha_{t}(l)& \mbox{otherwise}\
                \end{array}
\right.=\]
\[ = \left\{\begin{array}{ll}
                \alpha_r(k) & \mbox{if   }\ \
                   l=pqk,\ \mbox{for some } k\in \pf(r) ;\\
                %&\\
                \alpha_{s}(v)& \mbox{if   }\ \
                   l=pv,\ \mbox{for some } v\in \pf(s), \ q\not\prec v;\\
                \alpha_t(l) & \mbox{otherwise.}
                \end{array}
\right.\]

On the other side  we obtain $\la t,\alpha_t\ra(p;\la
s,\alpha_s\ra)=\la t(p;s),\gamma\ra$, where
\[ \gamma(m)=\left\{\begin{array}{ll}
                \alpha_s(v) & \mbox{if   }\ \
                   m=pv,\ \mbox{for some } v\in \pf(s) ;\\
                %&\\
                \alpha_{t}(m)& \mbox{otherwise}\
                \end{array}
\right.\] and $\la t(p;s),\gamma\ra=\la t(p;s)(pq;r),\delta\ra$,
where
\[ \delta(l)=\left\{\begin{array}{ll}
                \alpha_r(k) & \mbox{if   }\ \
                   l=pqk,\ \mbox{for some } k\in \pf(r) ;\\
                %&\\
                \gamma(l)& \mbox{otherwise}\
                \end{array}
\right.=\]
\[ = \left\{\begin{array}{ll}
                \alpha_r(k) & \mbox{if   }\ \
                   l=pqk,\ \mbox{for some } k\in \pf(r) ;\\
                %&\\
                \alpha_{s}(v)& \mbox{if   }\ \
                   l=pv,\ \mbox{for some } v\in \pf(s), \ q\not\prec v;\\
                \alpha_t(l) & \mbox{otherwise.}
                \end{array}
\right.\] Clearly, $\delta=\beta$ and  $t(p;s(q;r))= t(p;s)(pq;r)$.
\end{proof}

The following example illustrates the positions, subterms and
positional composition of colored terms.
\begin{example} \label{e-colorAlgproduct}
Let $\tau=(2)$, $\F=\{f\}$. The colorations of
 terms in the example are presented as bold
superscripts of the operation symbols.

Let  $\la t,\at\ra=f^{\bf 1}(f^{\bf 1}(x_1,x_2),f^{\bf
2}(x_1,x_2))$, $ \la s,\alpha_s\ra=f^{\bf 3}(f^{\bf
2}(x_1,x_2),x_2)$ and $ \la r,\alpha_r\ra=f^{\bf 3}(x_1,x_2)$ be
three  colored terms of type $\tau$.
Then we have $Pos(t)=\{\varepsilon,1,2,11,12,21,22\},$ $\la sub_t(2),\at[2]\ra=f^{\bf 2}(x_1,x_2),$ $
    \la sub_t(12),\at[12]\ra=x_2$ and  $Sub_c(\la
 s,\alpha_s\ra)=\{\la
 s,\alpha_s\ra, \la
 f(x_1,x_2),\alpha_s[1]\ra, x_1, x_2\}.$

  For the positional composition we have\\
\begin{tabular}{ll}
& $\la t,\at\ra(2;\la s,\alpha_s\ra(12;\la r,\alpha_r\ra))=$\\
& $f^{\bf 1}(f^{\bf 1}(x_1,x_2),f^{\bf 3}(f^{\bf
2}(x_1,x_2),x_2))(212;\la r,\alpha_r\ra)=$\\
 &$f^{\bf 1}(f^{\bf 1}(x_1,x_2),f^{\bf 3}(f^{\bf 2}(x_1,f^{\bf
3}(x_1,x_2)),x_2))=$\\
& $\la t,\at\ra(2;\la s,\alpha_s\ra)(212;\la r,\alpha_r\ra).$
\end{tabular}

\end{example}

\section{Multi-hypersubstitutions and deduction of identities}\label{sec4}

\begin{definition}\label{d-multihyp}\ \cite{dks1}\
 Let $M$ be a submonoid of ${H}yp(\tau)$ and let $\rho$ be
  a mapping of $\nn$ into $M$ i.e.
$\rho:\nn\to M.$ Any such mapping is called a {\it
multi-hypersubstitution} of type $\tau$ over $M$.
\end{definition}

  Denote by $\sigma_q$ the image of $q\in \nn$ under
$\rho$ i.e. $\rho(q)=\sigma_q\in Hyp(\tau).$

Let $\sigma\in M$ and $\rho\in Mhyp(\tau,M)$. If there is a natural
number $q\in\nn$ with $\rho(q)=\sigma$, then we will write
$\sigma\in\rho$.

We will  define the extension $\overline\rho_{\at}$ of a
multi-hypersubstitution $\rho$ to the set  of colored subterms of a
term.
%%%%%%%%%%%%%%%%%%%%%%%%%%%%%%%%%%%%%%%%%%%%%%

Let $\la t,\at\ra$ be a colored term of type $\tau$, $\at\in C_t$,
with $p\in Pos(t)$,  $s=sub_t(p)$ and $\at(p)=m$. Then we set:
\\
\begin{tabular}{ll}
  &  $(i)$\   if $s=x_j\in X$, then
$\overline\rho_{\at}[s]:=x_j;$
\\
%\hspace*{.2cm}
&  $(ii)$ \  if $s=f(s_1,\ldots,s_n)$, then
$\overline\rho_{\at}[s]:=\sigma_{m}(f)
(\overline\rho_{\at}[s_1],\ldots,\overline\rho_{\at}[s_n]).$
\end{tabular}
\\
The extension of $\rho$ on $\alpha_s$
  assigns inductively
  a coloration $\overline\rho_{t}[\alpha_s]$ to the term $\overline\rho_{\at}[s]$,
   as follows:
\\
\begin{tabular}{l}
 $(i)$\ if $s=f(x_1,\ldots,x_n)$, then
$\overline\rho_{t}[\alpha_s](q):=m$ for all
$q\in\pf(\overline\rho_{\at}[s]);$
\\
 $(ii)$\  if $s=f(s_1,\ldots,s_n),\ \mbox{and }\
q\in\pf(\overline\rho_{\at}[s]),$ then
\end{tabular}
\[\overline\rho_t[\alpha_s](q) =\left\{\begin{array}{ll}
               m & \mbox{  if  }
                 \ \ q\in \pf(\sigma_m(f));\\ \vspace{-.3cm}
                 &~~~~~~~~\\
                 \overline\rho_t[{\alpha_{s_j}}](k)& \mbox{if }\  \
                q=l k,\ \mbox{for some }j,\  j\leq n,\\
                &
                 k\in \pf(\overline\rho_{\alpha_t}[s_j]), \  l\in
                 Pos^X(\sigma_m(f)).
\end{array}
\right.\]

\begin{definition} \label{d-extension}
 The {\it mapping} $\overline\rho$ on the set $ W^c_\tau(X)$  is defined as follows:

$(i)$\quad $\overline\rho[x]:=x$ for all $x\in X;$

$(ii)$\quad if $t=f(t_1,\ldots,t_n)$ and $\at\in C_t$, then
$\overline\rho[\la t,\alpha_t\ra]:=
\la\overline\rho_{\at}[t],\overline\rho_{t}[\at]\ra.$
\end{definition}

\begin{example} \label{e-colorMulti}
Let  $\la t,\at\ra$  and $ \la s,\alpha_s\ra$ be the  colored terms
of type $\tau$ from Example \ref{e-colorAlgproduct}.

Let  $\sigma_1(f)=f(x_2,x_1)$, $\sigma_2(f)=f(f(x_2,x_1),x_2)$ and
$\sigma_3(f)=f(x_1,x_2)$ be hypersubstitutions of type $\tau$ and
$M=\{\sigma_1,\sigma_2,\sigma_3,\ldots\}$ be a submonoid of
$Hyp(\tau)$. Let $\rho\in Mhyp(\tau,M)$ with $\rho(m):=\sigma_m.$
Then we obtain
\\
\begin{tabular}{ll}
& $\overline\rho[\la t,\alpha_t\ra]$\\ & $=\overline\rho[f^{\bf
1}(f^{\bf 1}(x_1,x_2),f^{\bf 2}(x_1,x_2))]$\\ & $=f^{\bf 1} (f^{\bf
2}(f^{\bf 2}(x_2,x_1),x_2),f^{\bf 1}(x_2,x_1))$
\end{tabular}\\
and $\overline\rho[\la s,\alpha_s\ra]=\overline\rho[f^{\bf 3}(f^{\bf
2}(x_1,x_2),x_2)]=f^{\bf 3}(f^{\bf 2}(f^{\bf 2}(x_2,x_1),x_2),x_2)$.
The image of $\la t,\alpha_t\ra$ under $\rho$ is shown in Figure
\ref{f-12}.

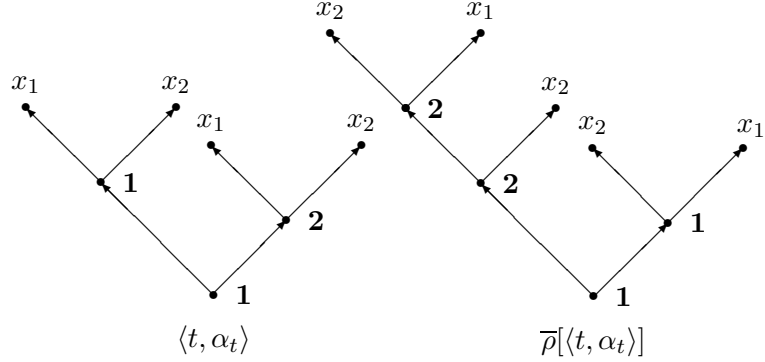
\begin{figure}
 \centering
%TeXCAD Picture [EX22.PIC]. Options:
%\grade{\overlinen}
%\emlines{\overlineff}
%\epic{\overlineff}
%\beziermacro{\overlinen}
%\reduce{\overlinen}
%\snapping{\overlineff}
%\quality{8.00}
%\graddiff{0.01}
%\snapasp{1}
%\zoom{9.5137}
\unitlength 1mm % = 2.85pt
\linethickness{0.4pt}
\ifx\plotpoint\undefined\newsavebox{\plotpoint}\fi % GNUPLOT compatibility
\begin{picture}(105.17,44.84)(0,0)
\put(33.67,7){\vector(-1,1){15}} \put(43.31,16.95){\vector(1,1){10}}
\put(43.31,16.95){\vector(-1,1){10}} \put(33.31,26.95){\circle*{1}}
\put(53.31,26.95){\circle*{1}} \put(43.31,16.95){\circle*{1}}
\put(33.67,7){\circle*{1}} \put(18.67,22){\vector(1,1){10}}
\put(33.39,6.97){\vector(1,1){10}} \put(18.67,22){\vector(-1,1){10}}
\put(8.67,32){\circle*{1}} \put(28.67,32){\circle*{1}}
\put(18.67,22){\circle*{1}} \put(37.67,7){\makebox(0,0)[cc]{{\bf
1}}} \put(22.67,22){\makebox(0,0)[cc]{{\bf 1}}}
\put(8.67,35){\makebox(0,0)[cc]{$x_1$}}
\put(28.67,35){\makebox(0,0)[cc]{$x_2$}}
\put(33.31,29.95){\makebox(0,0)[cc]{$x_1$}}
\put(53.31,29.95){\makebox(0,0)[cc]{$x_2$}}
\put(47.31,16.95){\makebox(0,0)[cc]{{\bf 2}}}
\put(33.67,1){\makebox(0,0)[cc]{$\langle t,\alpha_t\rangle$}}
\put(84.17,6.84){\vector(-1,1){15}}
\put(94.05,16.59){\vector(1,1){10}}
\put(94.05,16.59){\vector(-1,1){10}} \put(84.05,26.59){\circle*{1}}
\put(104.05,26.59){\circle*{1}} \put(94.05,16.59){\circle*{1}}
\put(84.17,6.84){\circle*{1}} \put(69.17,21.84){\vector(1,1){10}}
\put(84.17,6.71){\vector(1,1){10}}
\put(69.17,21.84){\vector(-1,1){10}} \put(59.17,31.84){\circle*{1}}
\put(79.17,31.84){\circle*{1}} \put(69.17,21.84){\circle*{1}}
\put(88.17,6.84){\makebox(0,0)[cc]{{\bf 1}}}
\put(73.17,21.84){\makebox(0,0)[cc]{{\bf 2}}}
\put(79.17,34.84){\makebox(0,0)[cc]{$x_2$}}
\put(84.05,29.59){\makebox(0,0)[cc]{$x_2$}}
\put(98.05,16.59){\makebox(0,0)[cc]{{\bf 1}}}
\put(84.17,.84){\makebox(0,0)[cc]{$\overline\rho[\langle
t,\alpha_t\rangle]$}} \put(59.17,31.84){\vector(1,1){10}}
\put(59.17,31.84){\vector(-1,1){10}} \put(49.17,41.84){\circle*{1}}
\put(69.17,41.84){\circle*{1}} \put(59.17,31.84){\circle*{1}}
\put(49.17,44.84){\makebox(0,0)[cc]{$x_2$}}
\put(69.17,44.84){\makebox(0,0)[cc]{$x_1$}}
\put(105.17,29.22){\makebox(0,0)[cc]{$x_1$}}
\put(63.17,31.84){\makebox(0,0)[cc]{{\bf 2}}}
\end{picture}

  \caption{Multi-hypersubstitution
  of colored terms}\label{f-12}
\end{figure}
\end{example}

\begin{proposition} \label{p-compTermMulti1}
Let $\rho\in Mhyp(\tau,M)$ be a multi-hypersubstitution and $t,s\in
W_\tau(X_n)$ with $\alpha_t\in C_t$, $\alpha_s\in C_s$ and $p\in
Pos(t)$ . Then
$$
\overline\rho[\la t,\alpha_t\ra(p;\la s,\alpha_s\ra)]=
\overline\rho[\la t,\alpha_t\ra(p;x_{j})](x_{j}\laa\overline\rho[\la
s,\alpha_s\ra]),
$$
for each $j,\ j>n$.
\end{proposition}
 \begin{proof}  We will use induction on the length $l(p)$ of
 the position
$p$.

First, let us observe that the
 case $l(p)=0$ is trivial.
So, our basis of induction is   $l(p)=1$. Then
$t=f(t_1,\ldots,t_{p-1},t_p,t_{p+1},\ldots,t_m)$ for some
$f\in\F_m$. Hence, for each $j,\ j>n$, we have
\\
\begin{tabular}{ll}
& $\overline\rho[\la t,\alpha_t\ra(p;\la
s,\alpha_s\ra)]=\overline\rho[\la
f(t_1,\ldots,t_{p-1},s,t_{p+1},\ldots,t_m),\alpha_{t'}\ra]=$\\
&$=\overline\rho[\la
f(t_1,\ldots,t_{p-1},x_{j},t_{p+1},\ldots,t_m)(x_{j}\laa s)
,\alpha_{t'}\ra]=$\\
& $=\overline\rho[\la
t,\alpha_t\ra(p;x_{j})](x_{j}\laa\overline\rho[\la s,\alpha_s\ra]),$
\end{tabular}
\\
where $\alpha_{t'}$ is a coloration of the term
$t'=f(t_1,\ldots,t_{p-1},s,t_{p+1},\ldots,t_m)$ for which
$\alpha_{t'}(q)=\alpha_t(q)$ when $q\in\pf(t)\setminus\{p\}$ and
$\alpha_{t'}(pq)=\alpha_s(q)$ when $q\in\pf(s).$

 Our inductive  supposition is
that when $l(p)<k$, the proposition is true, for some $k\in\nn$.

 Let $l(p)=k$ and $p=qi$ where
$q\in\nn_\F^*$ and $i\in\nn$. Hence $q\in\pf(t)$ and $l(q)<k$. Let
$u,v\in Sub(t)$ be subterms of $t$ for which $u=sub_t(p)$ and
$v=sub_t(q)$. Then we have
$v=g(v_1,\ldots,v_{i-1},v_i,v_{i+1},\ldots,v_j)$ with $v_i=u$ for
some $g\in\F_j.$ By the inductive supposition, for every $j,\ j>n$
and $l,\ l>n$, we obtain
\\
\begin{tabular}{ll}
& $\overline\rho[\la t,\alpha_t\ra(p;\la s,\alpha_s\ra)]=$
\\ & $=\overline\rho[\la t,\alpha_t\ra(q;x_{j})(x_{j}\laa \overline\rho[\la
v,\alpha_t[q]\ra(i;\la s,\alpha_s\ra)]=$
\\ & $=\overline\rho[\la
t,\alpha_t\ra(qi;x_{l})](x_{l}\laa\overline\rho[\la s,\alpha_s\ra])$
\\ & $=\overline\rho[\la
t,\alpha_t\ra(p;x_{j})](x_{j}\laa\overline\rho[\la s,\alpha_s\ra]).$
\end{tabular}

\end{proof}

A binary operation is defined in the set  $Mhyp(\tau,M)$  as
 follows:
\begin{definition} \label{d-prodMltiHyp}
Let  $\rho^{(1)}$ and $\rho^{(2)}$ be two multi-hypersubstitutions
over the submonoid $M$. Then the composition
$\rho^{(1)}\circ_{ch}\rho^{(2)}:\nn\to M$ maps each color $q\in \nn$
as follows:
$$(\rho^{(1)}\circ_{ch}\rho^{(2)})(q):=\rho^{(1)}(q)\circ_h\rho^{(2)}
(q):=\sigma^{(1)}_q\circ_h\sigma^{(2)}_q.$$
\end{definition}
\begin{lemma} \label{l-compMulti-hyp}
For every two  multi-hypersubstitutions $\rho^{(1)}$ and
$\rho^{(2)}$ over $M$ and for each colored term $\la t,\at\ra$ of
type $\tau$, it holds
$$\overline{(\rho^{(1)}\circ_{ch}\rho^{(2)})}[\la
t,\at\ra]=\overline\rho^{(1)}[\overline\rho^{(2)}[\la t,\at\ra]].$$
\end{lemma}

So, $ Mhyp(\tau,M)$ is a monoid, where
$\rho_{id}=(\sigma_{id},\sigma_{id},\ldots)$ is the identity
multi-hypersubstitution.

Let   $\A=\la A,\F\ra$ be an algebra of type $\tau.$ Each colored
term $\la t,\alpha_t\ra$ defines a term-operation
 on the set $A$, as follows $\la t,\at\ra^\A=t^\A$.

Let $\rho=(\sigma_1,\sigma_2,\ldots)$ be any multi-hypersubstitution
over a monoid $M$ of hypersubstitutions of type $\tau$, and let
$$\rho(\F)=\{t\in W_\tau(X)\ |\ \exists \ \sigma\in\rho\mbox{ and
}f\in\F\mbox{   such that   }t=\sigma(f)\}.$$

 The  algebra  $\rho[\A]:=\la A,\rho(\F)^{\rho[\A]}\ra$ is called a {\it
derived algebra} under the multi-hypersubstitution $\rho.$

 Let
$\R$ be a class of algebras of type $\tau$. The operator $\psi_M^c$
is defined as follows:\\
 \centerline{$\psi_M^c[\A]:=\{\rho[\A]\ |\ \rho\in Mhyp(\tau,M)\} \ \mbox{and}\
 \psi _{M}^{c}[\R]:=\{\psi_M^c[\A]\ |\ {\cal A}\in \R\}.$}
\begin{lemma}\label{l-ColAlg2}
For each $\la t,\alpha_t\ra\in W^c_\tau(X)$
it holds
$$\la t,\alpha_t\ra^{\rho[\A]}=\overline\rho[\la
t,\alpha_t\ra]^{\A}.$$
\end{lemma}
\begin{proof}  If  $t=x_j\in X$, then $\overline\rho[\la
t,\alpha_t\ra]^{\A}=x_j^{\A}\ \mbox{ and}\
 \la t,\alpha_t\ra^{\rho[{\A}]}=x_j^{\rho[{\A}]}=x_j^{\A}.$
Let us assume that $t=f(t_1,\ldots,t_n)$,
$\alpha_t(\varepsilon)=q\in \nn$, $\la f,q\ra^{\rho[{\A}]}= \rho(\la
f,q\ra)^\A$ and $\la
t_i,\alpha_{t_i}\ra^{\rho[\A]}=\overline\rho[\la
t_i,\alpha_{t_i}\ra]^{\A} $  for all  $i=1,\ldots,n$, then we have
$$\la t,\alpha_t\ra^{\rho[\A]}=\la f,q\ra^{\rho[\A]}(\la t_1,
\alpha_{t}[1]\ra^{\rho[\A]}
 ,\ldots,\la t_n,\alpha_{t}[n]\ra^{\rho[\A]})=$$
$$=\rho(\la f,q\ra)^\A(\overline\rho[\la t_1,
\alpha_{t}[1]\ra]^{\A}
 ,\ldots,\overline\rho[\la t_n,\alpha_{t}[n]\ra]^{ \A})=  \overline\rho[\la
t,\alpha_t\ra]^{\A}.$$ \end{proof}

 Let  $\rho$ be a multi-hypersubstitution. By  $\rho[ t\approx s]$
 we will denote
$\rho[ t\approx
s]:=\{\overline\rho_{\alpha_t}[t]\approx\overline\rho_{\alpha_s}[s]
\ |\ \alpha_t\in C_t, \ \alpha_s\in C_s\}.$ Let $\Sigma\subseteq
Id(\tau)$ be a set of identities of type $\tau$. The operator
$\chi_M^c$, is defined as follows:
\\
\begin{tabular}{ll}
 & $\chi_M^c[t\approx s]:=\{\rho[t\approx s]\ |\ \rho\in
Mhyp(\tau,M)\}$\\ & $ \mbox{and}\
\chi^c_M[\Sigma]:=\{\chi_M^c[t\approx s] \ |\ t\approx
s\in\Sigma\}.$
\end{tabular}
\begin{definition}\label{d-sch1}
An identity $t\approx  s\in Id\A$ in the algebra $\A$ is called an
{\it $M$-multi-hyperidentity } in $\A$, if for each
multi-hypersubstitution\linebreak
 $\rho\in Mhyp(\tau,M)$ and for every two colorations
 $\alpha_t\in C_t, \ \alpha_s\in C_s$, the identity
$\overline\rho_{\alpha_t}[t]\approx\overline\rho_{\alpha_s}[s]$ is
satisfied in $\A$.  When $t\approx s$ is an $M$-multi-hyperidentity
in $\A$ we will write $\A\ \models_{Mh}\  t\approx s,$  and the set
of all $M$-multi-hyperidentities in
 $\A$ is denoted by  $HC_MId\A.$
\end{definition}
 Algebras in which all identities are
$M$-multi-hyperidentities are called  {\it $M$-multi-solid
 } i.e.
an algebra $\A$ of type $\tau$ is  $M$-multi-solid, if
$\chi^c_{M}[Id\A]\subseteq Id\A.$  So, if  $V\subseteq\A lg(\tau)$
is a variety of type  $\tau$, it is called
 {\it  $M$-multi-solid}, when
 $\chi^c_{M}[IdV]\subseteq IdV.$
\begin{theorem}\label{t-ColorOperator11}
Let $\Sigma\subseteq Id(\tau)$ and $\R\subseteq \A lg(\tau)$. Then
\begin{enumerate}
\item[$(i)$]
$\chi_M[\Sigma ]\subseteq \chi_{M}^{c}[\Sigma ]$ \ for all $%
\Sigma \subseteq Id(\tau);$
\item[$(ii)$] $\psi_{M}[\R]=\psi_M^{c}[\R]$
\ for all  $\R\subseteq Alg(\tau).$
\end{enumerate}
\end{theorem}
\begin{proof}

$(i)$  Let $\sigma \in M$. Then we consider the
multi-hypersubstitution $\rho \in Mhyp(\tau,M )$ with
$\rho(q)=\sigma $
 for all $q \in \nn$. If \ $t\in W_{\tau }(X)$\ then we
will show that $\widehat {\sigma}[t]=\overline{\rho}_{\alpha_t}[t]$
for $\alpha_t\in C_t$. We will prove more,  that $\widehat
{\sigma}[r]:=\overline{\rho}_{\alpha_r}[r]$ for each $r \in Sub(t)$.
That will be proved  by induction on the depth of the term $r$.

 If $r\in X$ then
$\widehat{\sigma}[r]=\overline{\rho}_{\alpha_t}[r]=r$.

Let us assume that \ $r=f(r_{1},\ldots,r_{n})$
with $r_{1},\ldots,r_{n }%
\in Sub(t)$. Our inductive supposition is that
$\widehat{\sigma}[r_{k}]=\overline{\rho}_{\alpha_t}[r_{k}]$ for $1
\leq
k\leq n$. Then  we have\\
 \begin{tabular}{cc}
&$\overline{\rho}_{\alpha_t}[r]= \rho
(\at(\varepsilon))(f)(\overline{\rho}_{\alpha_t}[r_{1}],
\ldots,\overline{\rho}_{\alpha_t}[r_{n}])=\sigma (f)(\widehat{\sigma
}[r_{1}], \ldots,\widehat{ \sigma}[r_{n}])=$
\\
& $= \widehat{\sigma}[f(r_{1},\ldots,r_{n})]=\hat\sigma[r].$
\end{tabular}

 Let $t\approx s\in\Sigma$. Thus we have
$\hat\sigma[t]=\overline\rho_{\alpha_t}[t]$,
$\hat\sigma[s]=\overline\rho_{\alpha_s}[s]$ and
$$\hat\sigma[t]\approx \hat\sigma[s]\in \chi_M[\Sigma]\iff
\overline\rho_{\alpha_t}[t]\approx \overline\rho_{\alpha_s}[s]\in
\chi_M^c[\Sigma].$$ Hence $\hat\sigma[t]\approx \hat\sigma[s]\in
\chi_M^c[\Sigma].$

 $(ii)$ Let ${\cal A}\in \R$ and  $q\in\nn$ be the color
of $f$ in the fundamental colored term $\la
f(x_{1},\ldots,x_{n}),q\ra$. Let $\rho \in Mhyp(\tau,M )$. Then we
consider the hypersubstitution $\sigma \in M$ with $\sigma
(f)=\rho(q)(f)$ and from Lemma \ref{l-ColAlg2} we obtain
\[\la f,q\ra^{\rho \lbrack
{\cal A}]}={\overline\rho}[\la f(x_{1},\ldots,x_{n}),q\ra]^{\cal A}=
(\rho(q)(f)(x_{1},\ldots,x_{n}))^{\cal A}\]\[=\sigma
(f)(x_{1},\ldots,x_{n})^{\cal A}.\] This shows that $\psi _{M}^{c}[
{\cal A}]\subseteq \psi_{M}[{\cal A}].$

Let $\sigma \in Hyp(\tau)$,  $\A\in\R$ and $\rho \in Mhyp(\tau,M )$
with $\rho (q)=\sigma $ for all $q\in \nn$. From Lemma
\ref{l-ColAlg2} we obtain
$$\sigma(f)^{\cal A}=\hat\sigma[f(x_{1},\ldots,x_{n})]^{\cal
A}=\overline\rho[\la f(x_1,\ldots,x_n),q\ra]^{\cal A}=$$
$$=\la f(x_1,\ldots,x_n),q\ra]^{\rho [ {\cal A}]}=
\la f,q\ra^{\rho [ {\cal
A}]}.$$
 This shows that $\psi_{M}[{\cal A}]\subseteq \psi
_{M}^{c}[{\cal A}].$ Altogether we have $\psi _{M}^{c}[{\cal
A}]=\psi_{M}[{\cal A}]$ and thus $\psi_{M}^{c}[\R]=\psi_{M}[\R].$
\end{proof}

\begin{theorem}\label{p-rectbands1}
Each $M$-multi-solid variety is $M$-solid, but the converse
assumption is in general not true i.e. there are $M-$solid varieties
which are not $M-$multi-solid.
\end{theorem}
\begin{proof} Let $V$ be an $M$-multi-solid variety.
If $t\approx s\in\chi_M[IdV]$, then $t\approx s\in\chi_M^c[IdV]$
because of $(i)$ of Theorem \ref{t-ColorOperator11}. Hence every
$M$-multi-solid variety is $M$-solid one.

 Let $RB$ be the variety of rectangular bands, which is of
type   $\tau=(2)$. The identities satisfied in $RB$ are:
\\ \centerline{$Id_{RB}=$}
\\ \centerline{$\{f(x_1,f(x_2,x_3))\approx f(f(x_1,x_2),x_3) \approx
f(x_1,x_3),\  f(x_1,x_1)\approx x_1\}.$}

 In \cite{den51} it was
proved that $RB$ is a solid variety, which means that for each
identity $t\approx s$ with $Id_{RB}\models t\approx s$ and for each
hypersubstitution $\sigma$ we have $Id_{RB}\models
\hat\sigma[t]\approx \hat\sigma[s]$.

Let $M$, $\rho$ ,  $\la t,\at\ra$ and  $\la s,\alpha_s\ra$ be as in
Example \ref{e-colorMulti} i.e.
$$\la t,\at\ra=f^{\bf 1}(f^{\bf 1}(x_1,x_2),f^{\bf 2}(x_1,x_2))
\quad\mbox{ and}\quad \la s,\alpha_s\ra=f^{\bf 3}(f^{\bf
2}(x_1,x_2),x_2).$$ Since \\
\centerline{$t=f(f(x_1,x_2),f(x_1,x_2))\approx f(x_1,x_2)\approx
f(x_1,f(x_2,x_2))\approx f(f(x_1,x_2),x_2),$} it follows that
$Id_{RB}\models t\approx s$. Following the results in Example
\ref{e-colorMulti} we have
\\
\begin{tabular}{ll}
& $\overline\rho_{\alpha_t}[t]=f (f(f(x_2,x_1),x_2),f(x_2,x_1))$ and\\
& $\overline\rho_{\alpha_s}[s]=f (f(f(x_2,x_1),x_2),x_2).$
\end{tabular}\\
Thus we obtain $Id_{RB}\models\overline\rho_{\alpha_t}[t]\approx
f(x_2,x_1)$ and $Id_{RB}\models\overline\rho_{\alpha_s}[s]\approx
x_2.$ Hence
$Id_{RB}\not\models\overline\rho_{\alpha_t}[t]\approx\overline\rho_{\alpha_s}[s],$
and $RB$ is not  $M$-multi-solid. \end{proof}

The operators $\chi_M^{c}$ and $\psi_M^{c}$ are   connected by the
condition that
\\ \centerline{$ \psi_M^{c}[{\cal A}]\quad\mbox{\it satisfies }\quad u \approx v \iff
 {\cal A} \quad\mbox{\it satisfies }\quad \chi_M^{c} [u \approx v].
 $}
  $\psi_{M}^{c}$ and $\chi
_{M}^{c}$ are additive and closure operators.

Let  $\R\subset \A lg(\tau)$ and  $\Sigma\subset Id(\tau)$. Then we
set: \\ \centerline{$H{C}_MMod \Sigma:=\{{\cal A} \in \A lg(\tau) \
|\ t\approx s\in\Sigma\implies \A\models_{Mh} t\approx s\};$}\\
\centerline{ $H{C}_MId \R:=\{  {t\approx s} \in Id(\tau)\ |\
\chi_M^c[t\approx s] \subseteq
 Id \R\};$} \\
 and
\\ \centerline{$HC_MVar\R:=HC_MModHC_MId\R.$}
A set $\Sigma$ of identities of type $\tau$ is called an {\it
$M-$multi-hyperequational theory } if $\Sigma=HC_MIdHC_MMod\Sigma$
and a class $\R$ of algebras of type $\tau$ is called an {\it
$M-$multi-hyperequational class } if $\R=HC_MModHC_MId\R$.

\begin{theorem}\label{t-baseColor}
Let $\Sigma\subseteq Id(\tau)$ and $\R\subseteq \A lg(\tau)$. If
$\Sigma=HC_MId\R$ and $\R=HC_MMod\Sigma$, then $\R=Mod\Sigma$ and
$\Sigma=Id\R.$
\end{theorem}
\begin{proof}\ Let $\Sigma=HC_MId\R$ and $\R=HC_MMod\Sigma$.
 First, we will
prove that $\Sigma=\chi_M^c[\Sigma]$ and $\R=\psi_M^c[\R]$.

$\chi_M^c$ and $\psi_M^c$ are closure operators  and hence we have
$\Sigma\subseteq\chi_M^c[\Sigma]$ and $\R\subseteq\psi_M^c[\R]$.

We will prove the converse inclusions.
 Let $r\approx
v\in\chi_M^c[\Sigma]$. Then there is an identity $t\approx
s\in\Sigma$ with $r\approx v\in\chi_M^c[t\approx s]$, i.e.
$\chi_M^c[r\approx v]\subseteq \chi_M^c[\chi_M^c[t\approx
s]]\subseteq \chi_M^c[t\approx s].$
 From $t\approx s\in\Sigma$ it follows
$t\approx s\in HC_MId\R$, $\chi_M^c[t\approx s]\subseteq Id\R$ and
$\chi_M^c[r\approx v]\subseteq \chi_M^c[t\approx s]\subseteq Id\R.$
 Hence $r\approx v\in HC_MId\R$ and $r\approx v\in\Sigma$, i.e.
$\Sigma=\chi_M^c[\Sigma].$

Let $\A\in\psi_M^c[\R]$. Then there is an algebra $\B\in\R$ with
$\A\in\psi_M^c[\B]$, i.e. $\A=\rho[\B]$ for some $\rho\in
Mhyp(\tau,M)$.
 Hence $\psi_M^c[\A]= \psi_M^c[\rho[\B]]\subseteq
\psi_M^c[\B].$
 Since $\B\in\R=HC_MMod\Sigma$ we have $\chi_M^c[\Sigma]\subseteq Id\B$ i.e.
 $\psi_M^c[\B]\subseteq Mod\Sigma$  and
$\psi_M^c[\A]\subseteq \psi_M^c[\B]\subseteq Mod\Sigma.$
 Hence $\A\in HC_MMod\Sigma$ and $\A\in\R$, i.e.
$\R=\psi_M^c[\R].$

Now, we obtain
\\
\centerline{
 $Mod\Sigma=\{\A\ |\ \Sigma\subseteq Id\A\}=\{\A\ |\
\chi_M^c[\Sigma]\subseteq Id\A \}= HC_MMod\Sigma=\R$\  and}\\
\centerline{ $Id\R=\{r\approx v\ |\ \R\models r\approx
v\}=\{r\approx v\ |\ \psi_M^c[\R]\models r\approx v \}=
HC_MId\R=\Sigma.$}  \\
Hence $\R=Mod\Sigma$ and $\Sigma=Id\R.$
 \end{proof}

\begin{proposition}\label{p-main2}
 Let $\A\in\A lg(\tau)$, $\R\subseteq \A lg(\tau)$,
  $t\approx s\in Id(\tau)$,
 $\Sigma\subseteq Id(\tau)$, $\alpha_t\in C_t$ and $\alpha_s\in C_s$ and let $\rho$ be a
multi-hypersubstitution over $M$. Then we have:
\\
\begin{tabular}{l}
$(i)$\ $\A\models \overline\rho_{\alpha_t}[t]\approx
\overline\rho_{\alpha_s}[s] \iff \rho[\A]\models
t\approx s$ ;\\
$(ii)$\  $\A\models_{Mh}t\approx s
 \iff \chi_M^c[
t\approx s]\subseteq Id\A \iff \psi^c_M[\A]\models t\approx s;$
\\ $(iii)$\  $\Sigma\subseteq HC_MId{\mathcal R}
 \iff \chi_M^c[\Sigma]\subseteq Id{\mathcal R}
  \iff \Sigma\subseteq Id\psi^c_M [{\mathcal R}]$
  %\\ &$\iff \chi_M^c[\Sigma]\subseteq HC_MId\R;$
\\ $(iv)$\  ${\mathcal R}\subseteq HC_MMod\Sigma \iff
\psi^c_M[{\mathcal R}]\subseteq Mod\Sigma
  \iff {\mathcal R}\subseteq Mod\chi_M^c[\Sigma]$
%  \\ &$\iff \psi_M^c[\R]\subseteq HC_MMod\Sigma;$
\\ $(v)$\  $ HC_MId{\mathcal R} \subseteq Id{\mathcal R};$
\\ $(vi)$\  $ Var{\mathcal R} \subseteq HC_MVar{\mathcal R};$
\\
(vii) \   The pair $\la HC_MId,HC_MMod\ra$ forms a Galois
connection.
\end{tabular}
\end{proposition}
\begin{proof}
$(i)$ follows from Theorem \ref{t-baseColor} and  Lemma
\ref{l-ColAlg2}.

$(ii)$ Let $\rho$ be an arbitrary multi-hypersubstitution over the
monoid $M$. Let us assume  that  $\A\models_{Mh}t\approx s.$ Hence
 $\A\models \overline\rho_{\alpha_t}[t]\approx\overline\rho_{\alpha_s}[s]$
    i.e.  $\A\models\chi_M^c[ t\approx s]$
   and
    $\psi^c_M[\A]\models t\approx s.$

$(iii)$ and $(v)$ follow from  $(ii).$

$(iv)$ follows from   $(iii).$

$(vi)$ From   $(ii)$, it follows that $HC_MVar{\mathcal R}$ is a
variety. Then  \\ \centerline{${\mathcal R}\subseteq HC_MVar{\mathcal R}$, and $
Var{\mathcal R} \subseteq HC_MVar{\mathcal R}.$}

$(vii)$  follows from  Theorem \ref{t-baseColor}.
 \end{proof}

\begin{proposition}\label{p-main}
For every  $\Sigma\subseteq Id(\tau)$ and $\R\subseteq \A
lg(\tau)$ the following hold:

\begin{tabular}{ll}
(i)\ &$\chi_M^c[HC_MId\R]=HC_MId\R=
   Id\psi_M^c[\R]$;\\
(ii)\  &$\psi_M^c[HC_MMod\Sigma]=HC_MMod\Sigma=
   Mod\chi_M^c[\Sigma]$;\\
 (iii)\ &
$Var\psi_M^c[\R]=HC_MVar\R=ModHC_MId\R$;\\
(iv)\ & $Mod\chi_M^c[\Sigma]=HC_MMod\Sigma=IdHC_MMod\Sigma;$
\\
(iv)\ & $Id\psi_M^c[\R]=HC_MId\R=ModHC_MId\R.$
\end{tabular}
\end{proposition}
\begin{proof}\
$(i)$\ From Proposition \ref{p-main2} we
obtain \\
\begin{tabular}{ll}
&$HC_MId\R\subseteq HC_MId\R$\\
& $\implies\chi_M^c[\chi_M^c[HC_MId\R]]\subseteq
\chi_M^c[HC_MId\R]\subseteq
Id\R $\\
& $ \implies\chi_M^c[HC_MId\R]\subseteq HC_MId\R.$
\end{tabular}
\\
The converse inclusion is obvious.

$(ii)$ can be proved in an analogous way as $(i)$.

$(iii)$\ We have consequently,
\\
\begin{tabular}{ll}
& $HC_MVar\R=HC_MModHC_MId\R$\\
& $= Mod\chi_M^c[HC_MId\R]= ModHC_MId\R$\\
& $=ModId\psi_M^c[\R]=Var\psi_M^c[\R].$
\end{tabular}

 $(iv)$ and $(v)$ can be proved in an analogous way as $(iii)$.
 \end{proof}

For given set $\Sigma$ of identities the set $IdMod\Sigma$ of all
identities satisfied in the variety $Mod\Sigma$ is the deductive
closure of $\Sigma$, which is the smallest fully invariant
congruence  containing $\Sigma$ (see \cite{bir1,bur,mck}).

A remarkable fact is,  that there exists a variety $V\subset \A
lg(\tau)$ with $IdV=\Sigma$ if and only if $\Sigma$ is  a fully
invariant congruence  \cite{bur}.

 A congruence $\Sigma\subseteq Id(\tau)$ is called a fully
 invariant congruence if it additionally satisfies  the
following axioms (some authors call them ``deductive rules",
``derivation rules", ``productions" etc.):
\\ \hspace*{1cm}
$(i)$\quad {(variable inductive substitution)} \\ $(t\approx s\in
\Sigma)\ \&\ (r\in W_\tau(X))\ \&\ (x\in var(t))\  \implies\ t(x\laa
r)\approx s(x\laa r)\in\Sigma$;
\\ \hspace*{1cm}
$(ii)$\quad {(term positional replacement)} \\
$(t\approx s\in \Sigma)\ \&\ (r\in W_\tau(X))\ \&\ (sub_r(p)=t)\
\implies\ r(p;s)\approx r\in\Sigma$.

For any set of identities $\Sigma$ the smallest fully
 invariant congruence containing $\Sigma$ is called the $D-$closure of $\Sigma$
and it is denoted by $D(\Sigma).$

In \cite{den51} totally invariant congruences are studied as fully
invariant congruences which preserve the hypersubstitution images
i.e.  if $t\approx s\in\Sigma$ then
$\hat\sigma[t]\approx\hat\sigma[s]\in\Sigma$ for all $\sigma\in
Hyp(\tau).$

We  extend that results, to the case of multi-hypersubstitutions
over a given submonoid $M\subseteq Hyp(\tau).$
\begin{definition}
\label{d-chi} A  fully
 invariant congruence $\Sigma$  is $Mh$-deductively closed
if it additionally satisfies
\\ \hspace*{1cm}
$Mh_1$\quad {\it (Multi-Hypersubstitution)}

\begin{tabular}{c} $(t\approx s\in \Sigma)\  \&\ (\rho\in Mhyp(\tau,M))\  \implies\ \rho[t\approx
s]\subseteq\Sigma$.
\end{tabular}
\end{definition}

For any set of identities $\Sigma$ the smallest $Mh-$deductively
closed set containing $\Sigma$ is called the $Mh-$closure of
$\Sigma$ and it is denoted by $Mh(\Sigma).$ It is clear that for
each fully invariant congruence $\Sigma$ we have
$Mh(\Sigma)=\chi_M^c[\Sigma].$

 Let $\Sigma$
be a set of identities of type $\tau.$ For $t\approx s\in Id(\tau)$
we say $\Sigma\vdash_{Mh} t\approx s $ (``$\Sigma$ $Mh$-proves
$t\approx s$") if there is a sequence of identities $t_1\approx
s_1,\ldots,t_n\approx s_n$, such that each identity belongs to
$\Sigma$ or is a result of applying any of the derivation rules of
fully invariant congruence
   or $Mh_1$-rule
  to previous identities
in the sequence and the last identity $t_n\approx s_n$ is
$t\approx s.$

%%%%%%%%%%%%%%%%%%%%%%%%%%%%%%%%%

Let $t\approx s$ be an identity and $\A$ be an algebra of type
$\tau$. Then $\A\models_{Mh} t\approx s$ means that $\A\models
Mh(t\approx s)$ (see Definition \ref{d-sch1}).

Let $\Sigma\subseteq Id(\tau)$ be a set of  identities and $\A$ be
an algebra of type $\tau$. Then $\A\models_{Mh} \Sigma$ means that
$\A\models Mh(\Sigma)$.
 For $t,s\in W_\tau(X)$ we say
$\Sigma\models_{Mh} t\approx s$ (read: ``$\Sigma$ $Mh-$yields
$t\approx s$") if, given any algebra $\B\in \A lg(\tau)$,
$$\B\models_{Mh} \Sigma\quad\Ra\quad \B\models_{Mh} t\approx s.$$
\begin{remark}\label{r-multi-solid}
In a more general case we have $\chi_M^c[\Sigma]\subseteq
Mh(\Sigma)$ and there are examples when $\chi_M^c[\Sigma]\neq
Mh(\Sigma)$.

It is easy to see that each $Mh-$deductively closed set is a totally
invariant congruence \cite{den51}. On the other side from Theorem
\ref{p-rectbands1} it follows that the totally invariant congruence
$Id_{RB}$ is not $Mh-$deductively closed.
\end{remark}

\begin{lemma}\label{l-VdashMod}
For any set $\Sigma\subseteq Id(\tau)$ of identities and $t\approx
s\in Id(\tau)$ the following equivalences hold:
$$\Sigma\vdash_{Mh} t\approx s \iff \chi_M^c[\Sigma]\vdash
t\approx s\iff Mh(\Sigma)\vdash t\approx s.$$
\end{lemma}
\begin{theorem}\label{t-complBirkhoff}
{\bf (Completeness Theorem for Multi-hyperequational Logic.) } For
$\Sigma\subseteq Id(\tau)$ and $t\approx s \in Id(\tau)$ we have:
$$\Sigma\models_{Mh}t\approx s\iff \Sigma\vdash_{Mh}t\approx s.$$
\end{theorem}
\begin{proof} \ From
$$HC_MIdHC_MMod\Sigma=IdHC_MMod\Sigma=IdMod\chi_M^c[\Sigma]$$ we
obtain that  $\Sigma\models_{Mh}t\approx s$ is equivalent to
 $t\approx
s\in HC_MIdHC_MMod\Sigma$.   From Theorem 14.19 \cite{bur} we have
$$HC_MIdHC_MMod\Sigma\models t\approx s\ \mbox{ and }\ 
HC_MIdHC_MMod\Sigma\vdash t\approx s.$$ Hence
$\chi_M^c[\Sigma]\vdash t\approx s$ and $\Sigma\vdash_{Mh}t\approx
s.$

The converse implication follows from the fact that $Mh(\Sigma)$ is
a fully invariant congruence which is  closed under the rule $Mh_1$.
\end{proof}

\begin{corollary}\label{c-completeness2}
Let $M$ be a submonoid of $Hyp(\tau)$. Then the class of all
$M-$multi-solid varieties of type $\tau$ is a complete sublattice
of the lattice $\mathcal L(\tau)$ of all varieties of type $\tau$
and dually, the class of all $M-$multi-hyperequational theories of
type $\tau$ is a complete sublattice of the lattice $\E(\tau)$ of
all equational theories (fully invariant congruences) of type
$\tau$.
\end{corollary}
\begin{lemma}\label{l-ModelsMod}
For any set $\Sigma\subseteq Id(\tau)$ of identities and $t\approx
s\in Id(\tau)$ the following equivalence holds:
$$\Sigma\models_{Mh} t\approx s \iff Mh(\Sigma)\models t\approx
s.$$
\end{lemma}
\begin{proof}

$``\Ra"$ Let $\Sigma\models_{Mh} t\approx s$. We have to prove
$Mh(\Sigma)\models t\approx s$. Let $\A\in\A lg(\tau)$ be an
algebra for which $\A\models Mh(\Sigma)$. This implies
$\A\models_{Mh}\Sigma$ and $\A\models_{Mh} t\approx s$, because of
$\Sigma\models_{Mh} t\approx s$. Hence $\A\models Mh(t\approx s)$.
On the other side, we have $t\approx s\in Mh(t\approx s)$ and
$\A\models t\approx s$.

$``\La"$ Let $Mh(\Sigma)\models t\approx s.$ Since
$\chi_M^c[\Sigma]\subseteq Mh(\Sigma)$ and from Proposition
\ref{p-main2} $(ii)$ we have $\Sigma\models_{Mh}t\approx s.$
\end{proof}
\begin{corollary}\label{c-ModelsMod}
For any set $\Sigma\subseteq Id(\tau)$ of identities and $t\approx
s\in Id(\tau)$ it holds: $$\Sigma\models_{Mh} t\approx s \iff
\chi_M^c[\Sigma]\models t\approx s.$$
\end{corollary}

\begin{example}\label{e-vain}
Let $\tau$ be an arbitrary type.
\par
$(i)$\ Let\   $t\in W_\tau(X)$  be a term of type $\tau.$ Let us
consider the following two functions $Left:W_\tau(X)\to X$ and
$Right:W_\tau(X)\to X$, which assign to each term  $t$ the
leftmost and the rightmost variable of  $t$. For instance, if
$t=g(f(x_1,x_2),x_3,x_4),$ then $ Left(t)=x_1$ and $
Right(t)=x_4.$

Let us denote by $K_1\subset Hyp(\tau)$ the set of all
hypersubstitutions which preserve the functions   $Left$ and
$Right$ i.e. $\sigma\in K_1$, iff for all  $t\in W_\tau(X)$ we
have $Left(t)=Left(\hat\sigma[t])\quad\mbox{ and}\quad
 Right(t)=Right(\hat\sigma[t]).$
Then $ K_1$ is a submonoid of $Hyp(\tau)$. Let us consider the
monoid   $Mhyp(\tau, K_1)$  of multi-hypersubstitutions, generated
by  $K_1$. It is a submonoid of $Mhyp(\tau, Hyp(\tau) ).$ The
variety  $RB$ of rectangular bands is  $K_1$-multi-solid.

$(ii)$\ Let $Vstr:W_\tau(X)\to X^*$ be a mapping, which assigns to
each term   $t$ the string of the variables in $t$. For instance,
if
$$t=f(x_1,g(f(x_1,x_2),x_3,x_2),x_4),\quad\mbox{ then}\quad
Vstr(t)=x_1x_1x_2x_3x_2x_4.$$ This mapping is defined inductively
as follows: if $t=x_j\in X$, then \\ $Vstr(t):=x_j;$ and
 if
$t=f(t_1,\ldots,t_n)$, then\\ $Vstr(t):=Vstr(t_1) Vstr(t_2) \ldots
Vstr(t_n).$

 Let  $K_2\subset Hyp(\tau)$ be the set of all
hypersubstitutions which preserve  $Vstr$ i.e. $\sigma\in K_2$, if
and only if for each  $t\in W_\tau(X)$ it holds
$Vstr(t)=Vstr(\hat\sigma[t]).$
 Then  $ K_2$ is a submonoid of  $Hyp(\tau)$.

Let us consider the monoid   $Mhyp(\tau, K_2)$ of
multi-hypersubstitutions, generated by  $K_2$. It is a submonoid
of $Mhyp(\tau, Hyp(\tau) )$ and
 $K_2\subseteq K_1$.

 It is not difficult to
prove that the variety $RB$
 is   $K_2$-multi-solid.

Finally, let us note that if $\Sigma $ is the set of identities
satisfied in $RB$, and if we add to $K_1$ and $K_2$ the
hypersubstitution $\sigma\in Hyp(\tau)$ with
$\sigma(f)=f(x_2,x_1)$, then $M_1:=K_1\cup\{\sigma\}$ and
$M_2:=K_2\cup\{\sigma\}$ are monoids, again such that
$\chi_{M_1}^c[\Sigma]=Id(\tau)$, but $\chi_{M_2}^c[\Sigma]\neq
Id(\tau)$.
\end{example}

\section{Tree automata realization}\label{sec7}~~~

We consider an automata realization of the multi-hypersubstitutions.
This concept will allow to use computer programmes in the case of
finite monoids of hypersubstitutions to obtain the images of terms
under multi- hypersubstitutions.

In Computer Science terms are used as data structures and
  they are called trees.  The operation symbols are labels of
the internal nodes of trees and variables are their leaves.

 The
concept of tree automata was introduced  in the 1960s in various
papers such as   \cite{Tha1}.    G\'ecseg \&   Steinby's book
\cite{ges} is a good survey of the theory of tree automata and
\cite{Com} is a development of this theory. Tree automata are
classified  as tree recognizers and tree transducers. Our aim  is to
define tree automata which interpret the application of
multi-hypersubstitutions of the colored terms of a given type.

\begin{definition}\label{d-treeRhoTrans}
A colored tree transducer of type $\tau$    is a tuple
$\underline{A}$ = $\la X, \F,
 P\ra$ where  as usual $X$ is a set of variables, $\F$ is a set
 of operation symbols and $P $ is a finite set of {\it
productions~(rules~ of~ derivation)} of the forms
\\
 \hspace{3mm}
\begin{tabular}{l}
(i)\  $x \rightarrow x,\ x\in X;$\\
(ii) \ $\la f,q\ra (\xi_1, \cdots , \xi_{n}) \rightarrow  \la
r,\alpha_r\ra (\xi_1,  \cdots, \xi_{n}),$ \\  with $\la
r,\alpha_r\ra (\xi_1, \cdots, \xi_{n}) \in W_{\tau}^c(X \cup
\chi_m)$, $\la f,q\ra \in \F^c,$ $\xi_1,
\cdots, \xi_{n} \in \chi_{n}$,\\
 where   $\chi_{n} = \{\xi_1,
\cdots, \xi_{n} \}$ is an auxiliary alphabet.
\end{tabular}
\\
(All auxiliary variables $\xi_j$ belong to a set $\chi_m =\{\xi_1,
\cdots, \xi_m \}$ where $m=maxar$ is the maximum of the arities of
all operation symbols in $\F$.)
\end{definition}

%%%%%%%%%%%%%%%%%%%%%%%%%
\begin{definition} Let $M\subseteq \H yp(\tau)$ be a submonoid of
$\H yp(\tau)$ and $\rho\in Mhyp(\tau,M)$ be a
multi-hypersubstitution over $M$. A colored tree transducer
$\underline{A}^\rho$ = $\la X, \F,
 P\ra$ is called {\it $Mh-$transducer} over $M$, if for the rules (ii) in $P$
we have $r=\sigma(q)(f)$ and $\alpha_r(p)=q$ for all $p\in \pf(r)$.
\end{definition}

 The $Mh-$transducer  $\u A^\rho$
 runs over a colored term $\la t,\at\ra$  starting  at the
leaves of $t$ and moves downwards, associating along the run a
resulting colored term (image)  with each subterm inductively: if
$t=x\in X$ then the $Mh-$transducer  $\u A^\rho$  associates with
$t$ the term $x\in W^c_\tau(X)$,  if $x \rightarrow x\in P;$ if $\la
t,\at\ra=\la f,q\ra(\la t_1,\alpha_{t_1}\ra\ldots,\la
t_n,\alpha_{t_n}\ra)$  then with $\la t,\at\ra$ the $Mh-$transducer
$\u A^\rho$ associates the colored term
 $\la s,\alpha_s\ra=\la
u,\alpha_{u}\ra(\la t_1,\alpha_{t_1}\ra\ldots,\la
t_n,\alpha_{t_n}\ra),$ if $ \la f,q\ra({\xi}_1, \dots, {\xi}_{n})
\rightarrow \la u,\alpha_u\ra({\xi}_1, \cdots, {\xi}_{n})\in P$,
where $u=\sigma(q)(f)$ and $\alpha_u(p)=q$ for all $p\in \pf(u).$

For trees $\la t,\alpha_t\ra,$ and $\la  s,\alpha_s\ra$ we say {\it
$\la t,\alpha_t\ra$ ~directly derives~ $\la  s,\alpha_s\ra$
 ~by ~$\underline A^\rho$},
if $\la  s,\alpha_s\ra$ can be obtained from $\la t,\alpha_t\ra$ by
 replacing of an occurrence of a
subtree  $\la f,q\ra(\la r_1,\alpha_{r_1}\ra, \cdots, \la
r_n,\alpha_{r_n}\ra),
 $ in $\la t,\alpha_t\ra$
 by \\
 $\la
u,\alpha_u\ra(\la r_1,\alpha_{r_1}\ra, \cdots, \la
r_n,\alpha_{r_n}\ra)\in W_{\tau}^c(X \cup \chi_m)$.

If $\la t,\alpha_t\ra$ directly derives $\la s,\alpha_s\ra$ in
$\underline A^\rho$, we write $\la t,\alpha_t\ra
\rightarrow_{\underline A^\rho} \la s,\alpha_s\ra$. Furthermore, we
say  {\it  $\la t,\alpha_t\ra$ ~derives $\la s,\alpha_s\ra$ ~in
~$\underline A^\rho$}, if there is a sequence $$\la
t,\alpha_{t}\ra\rightarrow_{\underline A^\rho} \la
s_1,\alpha_{s_1}\ra \rightarrow_{\underline A^\rho} \la
s_2,\alpha_{s_2}\ra \rightarrow_{\underline A^\rho} \cdots
\rightarrow_{\underline A^\rho} \la s_n,\alpha_{s_n}\ra = \la
s,\alpha_{s}\ra$$ of direct derivations  or if $\la t,\alpha_{t}\ra =
\la s,\alpha_{s}\ra$. In this case we write $\la t,\alpha_{t}\ra
 \Rightarrow^*_{\underline A^\rho}\la s,\alpha_{s}\ra $.
Clearly, $\Rightarrow^*_{\underline A^\rho}$ is the reflexive and
transitive closure of $\rightarrow_{\underline A^\rho}$.

Let us denote by $Thyp(\tau,M)$ the set of all $Mh-$transducers of
type $\tau$ over $M.$

A term $\la t,\at\ra$ is translated to the term $\la s,\alpha_s\ra$
by the $Mh-$transducer  $\u A^\rho$ if there exists a run of $\u
A^\rho$ such that it associates with $\la t,\at\ra$ the colored term
$\la s,\alpha_s\ra$. In this case we will write
 $ \u
A^\rho(\la t,\alpha_t\ra)=\la s,\alpha_s\ra.$
\begin{lemma}\label{lem}
$ \u A^\rho(\la t,\alpha_t\ra)=\la s,\alpha_s\ra\iff
\overline\rho[\la t,\alpha_t\ra]=\la s,\alpha_s\ra.$
\end{lemma}
\begin{proof}
For a variable-term $x$ we have $\u A^{\rho}(x)=x=\overline\rho[x]$.
Suppose that for $i=1,\ldots,n$ we have $\u A^{\rho}(\la
t_i,\alpha{t_i}\ra)=\overline\rho[\la t_i,\alpha{t_i}\ra]$. Then we
obtain

$$\u A^{\rho}(\la f,q\ra(\la t_1,\alpha{t_1}\ra,\ldots,\la
t_n,\alpha{t_n}\ra))=$$  $$\u A^{\rho}(\la f,q\ra)(\u A^{\rho}(\la
t_1,\alpha{t_1}\ra),\ldots,\u A^{\rho}(\la
t_n,\alpha{t_n}\ra))=\overline\rho[\la
t,\alpha_t\ra]=\overline\rho[\la t,\alpha_t\ra],$$ where
$\alpha_t(\varepsilon)=q$ and $\alpha_t[i]=\alpha_{t_i}$ for
$i=1,\ldots,n$.
\end{proof}
The product (superposition) of two $Mh-$transducers $\u A^{\rho_1}$
and $\u A^{\rho_2}$ is defined by the following equation
$$\u A^{\rho_1}\circ \u A^{\rho_2}(\la t,\alpha_t\ra):=\u
A^{\rho_1}(\u A^{\rho_2}(\la t,\alpha_t\ra)).$$

\begin{lemma}\label{l-TransCompos}
Let $M\subseteq \H yp(\tau)$ be a monoid of hypersubstitutions. Then
the superposition of $Mh-$transducers of a given type is associative
i.e.
$$\u A^{\rho_1}\circ(\u A^{\rho_2}\circ \u A^{\rho_3})(\la t,\at\ra)=
(\u A^{\rho_1}\circ \u A^{\rho_2})\circ \u A^{\rho_3}(\la
t,\at\ra)$$ for all ${\rho_1}, {\rho_2},  {\rho_3}\in MHyp(\tau,M)$
and for all $\la t,\at\ra\in W^c_\tau(X)$.
\end{lemma}

\begin{theorem}\label{t-TransCompos}
Let $M\subseteq \H yp(\tau)$ be a monoid of hypersubstitutions. Then
the set $ Thyp(\tau,M)$ is a monoid which is isomorphic to the
monoid of all multi-hypersubstitutions over $M$ i.e.
$$Thyp(\tau,M)\cong Mhyp(\tau,M).$$
\end{theorem}
\begin{proof}
We define a mapping $\varphi: Mhyp(\tau,M)\to Thyp(\tau,M)$ by
$\varphi(\rho):=\u A^\rho.$

To show that $\varphi$ is a homomorphism we will prove that $\u
A^{\rho_1}\circ\u A^{\rho_2}=\u A^{\rho_1\circ_{ch}\rho_2}$, so that
$\varphi(\rho_1)\circ\varphi(\rho_2)=\varphi(\rho_1\circ_{ch}\rho_2)$.

We have $$\u A^{\rho_1}\circ\u A^{\rho_2}(\la t,\alpha_t\ra)=\u
A^{\rho_1}(\u A^{\rho_2}(\la t,\alpha_t\ra))=\u
A^{\rho_1}(\overline\rho_2[\la
t,\alpha_t\ra])=\overline\rho_1[\overline\rho_2[\la
t,\alpha_t\ra]]$$ $$=\overline{\rho_1\circ_{ch}\rho_2}[\la
t,\alpha_t\ra]=A^{\rho_1\circ_{ch}\rho_2}(\la t,\alpha_t\ra).$$

To see that $\varphi$ is one-to-one, let $\u A^{\rho_1}=\u
A^{\rho_2}$. Then from Lemma \ref{lem} for all $\la t,\alpha_t\ra\in
W_\tau^c(X)$ we have $\overline\rho_1[\la
t,\alpha_t\ra]=\overline\rho_2[\la t,\alpha_t\ra]$. Hence for all
$f\in \F$ and $q\in\nn$ we have $\overline\rho_1[\la
f,q\ra]=\overline\rho_2[\la f,q\ra]$ and therefore $\rho_1=\rho_2$.
\end{proof}


\begin{thebibliography} {20}
\bibitem{bir1} G. Birkhoff, {\it Lattice theory,} (3rd ed.)
 Amer. Math. Soc., Providence, 1967
\bibitem{bur} S. Burris and  H. Sankappanavar,
 {\it A Course in Universal Algebra},
The millennium edition, 2000
\bibitem{Com} H. Comon, M. Dauchet,
R. Gilleron, F. Jacquemard, D. Lugiez, S. Tison, M. Tommasi,{\it
Tree Automata, Techniques and Applications,} 1999,
http://www.grappa.univ-lille3.fr/tata/
\bibitem{dks1} K. Denecke, J. Koppitz and Sl. Shtrakov,
 Multi-Hypersubstitutions and Coloured Solid Varieties, {\it J.
Algebra and Computation,} J. Algebra and Computation, Volume 16,
Number 4, August, 2006, pp.797-815.
\bibitem{den51} K. Denecke, D. Lau, R. P\"oschel and D. Schweigert,
 Hyperidentities, Hyperequational Classes and Clone
Congruences,{\it General Algebra 7, Verlag
H\"older-Pichler-Tempsky,} Wien 1991, Verlag B.G. Teubner
Stuttgart, pp.97-118
\bibitem{ges} F. G\'ecseg, M. Steinby, {\it Tree Automata},
 Akad\'emiai
Kiad\'o, Budapest 1984
\bibitem{grac2} E. Gracz\'ynska, {\it On connection between identities and hyperidentities},
Bull.Sect.Logic 17(1988),34-41.
\bibitem{gra2} G. Gratzer, {\it Universal Algebra}, D. van Nostrand
Co., Princetown, 1968.
\bibitem{mck} R. McKenzie, G. Mc Nulty and W. Taylor, {\it Algebras,
Lattices, Varieties,} Vol. I, Belmont, California 1987.
\bibitem{Tha1} J. W. Thatcher and J.B. Wright,
 Generalized finite automata, {\it Notices Amer. Math. Soc.},
{\bf 12.} (1965), abstract No. 65T-649,820.

\end{thebibliography}
\end{document}